\documentclass{article}
\usepackage{amsfonts}

\newtheorem{theorem}{Theorem}

\def\C{{\mathbb C}}
\def\R{{\mathbb R}}
\def\beq{\begin{equation}}
\def\eeq{\end{equation}}
\def\const{{\mathrm{const}}}
\def\res{{\mathrm{res}}}
\def\nm{{\mathrm{nm}}}
\def\Sing{{\mathrm{Sing}\,}}

\begin{document}

\title{Orthogonal curvilinear coordinate systems corresponding to 
singular spectral curves
\thanks{The work is supported by RFBR (grant no. 03-01-00403) and 
the program of basic researches of
RAS ``Mathematical methods in nonlinear dynamics''. The first author
(A.E.M.) was also supported by the grant of President of Russian
Federation (grant MK-1017.2004.1) and the second author (I.A.T.) 
was also supported by 
Max Planck Institute of Mathematics.}}
\author{A.E. Mironov \thanks{Institute of Mathematics, 
630090 Novosibirsk, Russia; mironov@math.nsc.ru}
\and I.A. Taimanov \thanks{Institute of Mathematics, 630090
Novosibirsk, Russia; taimanov@math.nsc.ru}}
\date{}
\maketitle

\section{Introduction}

In this paper we study the limiting case of 
the Krichever construction of orthogonal curvilinear coordinate systems when 
the spectral curve becomes singular.

Theory of orthogonal curvilinear coordinate systems was very popular 
among differential geometers in the 19th century and in the first half of
the 20th century (Dupin, Gauss, Lame, Bianchi, Darboux) and the classification
problem was basically solved at the beginning of the 20th century
(see the book by Darboux \cite{D} which summarizes this stage of the 
development of this theory). 
These coordinate systems are interesting due to search of systems solved
by the separation of variables method and to modern problems of 
theory of hydrodynamical type systems and topological field theory
(Dubrovin, Novikov, Tsarev, Krichever, see references in \cite{K,Z}).

In the problem of explicit constructing such systems 
a breakthrough was achieved by Zakharov \cite{Z} who 
by using the dressing method first applied the methods of integrable
systems to this problem. Onto the finite-gap integration method this
approach was extended by Krichever \cite{K}. Therewith the initial data for 
a construction of such a system consist of a Riemann surface, i.e., 
the spectral curve, which in \cite{K} is assumed to be nonsingular and
some other additional quantities related to it. We briefly recall the 
constructions by Zakharov and Krichever in \S 2.

In the case when the spectral curve 
becomes singular and reducible such that all its 
irreducible components are smooth rational complex curves 
the procedure of constructing orthogonal curvilinear coordinate 
systems is crucially simplified and reduces to simple computations 
with elementary functions 
(see \S 3). Therewith we show how such well-known coordinate 
systems as the polar coordinates on the plane, the cylindrical coordinates 
in the three-space, and the spherical coordinates in $\R^n$ with 
$n \geq 3$ fit in this scheme (we expose these constructions together with 
constructions of some other coordinate systems in \S 4).

We remark that the inverse problems with such spectral curves were already 
studied in relations to their applications (see, for instance, the paper 
\cite{T2} on surface theory and the paper \cite{TC} on  the Hitchin system). 
Although this case is very special this paper also shows that  explicit
solutions corresponding to it are important for applications.

\section{Methods of constructing orthogonal curvilinear coordinates}

\subsection{Preliminary facts}

A curvilinear coordinate system 
$u=(u^1,\dots,u^n)$ in the Euclidean $n$-space 
${\mathbb R}^n$ is called $n$-orthogonal if the metric in these coordinates
takes the form
$$
 ds^2=H_1^2(du^1)^2+\dots+H_n^2(du^n)^2.
$$
Therewith the functions $H_j=H_j(u)$ are called the Lame coefficients
and the condition that the curvature tensor vanishes takes the form
\beq \label{1}
\frac{\partial^2 H_i} 
{\partial u^j \partial u^k} =
\frac{1}{H_j}\frac{\partial H_j}{\partial u^k}\frac{\partial H_i}{\partial u^j}
+ \frac{1}{H_k}\frac{\partial H_k}{\partial u^j}
\frac{\partial H_i}{\partial u^k},
\eeq 
\beq 
\label{2}
\frac{\partial}{\partial u^j}
\left(\frac{1}{H_j} \frac{\partial H_i}{\partial u^j}\right)
+
\frac{\partial}{\partial u^i}
\left(\frac{1}{H_i} \frac{\partial H_j}{\partial u^i}\right)
+
\sum_{k\ne i \ne j}\frac{1}{H_k^2}\frac{\partial H_i}{\partial u^k}
\frac{\partial H_j}{\partial u^k} = 0.
\eeq 
There are $\frac{n(n-1)(n-2)}{2}$ and $\frac{n(n-1)}{2}$
equations in the systems (\ref{1}) and (\ref{2}) respectfully,
the equations (\ref{1}) are equivalent to the condition that
$ R_{ijik}=0, j
\ne k$, and the equations (\ref{2}) are equivalent to
$R_{ijij}=0$.
Other components of the curvature tensor $R_{ijkl}$ 
always vanish for a diagonal metric. Hence the system of equations
(\ref{1})--(\ref{2}) for the Lame coefficients is strongly overdetermined. 
A general solution to this system is parameterized by 
$\frac{n(n-1)}{2}$
functions of two variables.

The order of the system (\ref{1})--(\ref{2}) is minimized by introducing 
the rotation coefficients
\beq
\label{3}
 \beta_{ij}=\frac{\partial_{u^i}H_j}{H_i}.
\eeq
Then the equations (\ref{1}) and (\ref{2}) take the form
\beq
\label{4}
 \partial_{u^k}\beta_{ij}=\beta_{ik}\beta_{kj},
\eeq
\beq
\label{5}
 \partial_{u^i}\beta_{ij}+\partial_{u^j}\beta_{ij}+\sum_{k\ne i,l}
 \beta_{ki}\beta_{kj}=0.
\eeq
Given a solution $\beta_{ij}$ to these equations, the Lame coefficients
are found from (\ref{3}) as a solution to the Cauchy problem
$$
H_i(0,\dots,0,u^i,0,\dots,0)=h_i(u^i).
$$
Therewith such a solution depends on the initial data for this problem,
i.e., on $n$ functions $h_i$ of one variable.We remark that the compatibility 
of (\ref{3}) is equivalent to (\ref{4}).

The immersion problem, i.e., the determination of the Euclidean coordinates
$x=(x^1,\dots,x^n)$ as the functions of $u=(u^1,\dots,u^n)$, reduces to solving
the overdetermined system of second order linear equations 
\beq \label{6}
\frac{\partial^2 x^k}{\partial u^i \partial u^j} =
\sum_{l}\Gamma_{ij}^l\partial_{u^l}x^k,
\eeq
and in our case the Christoffel symbols have the form
$$
 \Gamma_{ij}^k=0,\ i\ne j\ne k; \ \
 \Gamma_{kj}^k=\frac{\partial_{u^j}H_k}{H_k}; \ \
 \Gamma_{ii}^k=-\frac{H_i\partial_{u^k}H_i}{(H_k)^2},\ k\ne i.
$$
By (\ref{1}) and (\ref{2}), the 
system of equations (\ref{6}) is compatible
and determines an $n$-orthogonal curvilinear coordinate system up to 
motions of $\R^n$.

\subsection{Zakharov's method \cite{Z}}

The abstract $n$ waves problem has the form
\beq
\label{7}
 \sum_{i,j,k}\varepsilon_{ijk}(I_j\partial_{u^i}QI_k-I_iQI_jQI_k)=0,
\eeq 
where $Q(u)$ is the unknown $(n\times n)$-matrix function, the matrices
$I_j=I_j(u^j)$ are pairwise commuting, 
$\varepsilon_{ijk}=1$ for $i>j>k$ and changes the sign after an odd 
permutation of indices.

We take for $I_j$ the diagonal matrices with the diagonal 
$(0,\dots,0,1,0,\dots,0)$ (i.e., the unit is on the $j$-th place).
The equations (\ref{7}) take the form
$$
 \partial_{u^j}Q_{ik}=Q_{ij}Q_{jk},
$$
i.e., coincide with (4). Let us consider an auxiliary function
$\widetilde{Q}=\tilde{Q}(u,s)$, 
with $s = u^{n+1}$ an additional variable and
$I_{n+1}$ the unit matrix, satisfying the equations
\beq
\label{8}
 \partial_{u^j}[I_i,\widetilde{Q}]-\partial_{u^i}[I_j,\widetilde{Q}]+
I_i\partial_{s}\widetilde{Q}I_j-
 I_j\partial_{s}\widetilde{Q}I_i-[[I_i,\widetilde{Q}],[I_j,\widetilde{Q}]]=0,
\eeq
where $i,j=1,\dots,n+1$. It is shown that if 
$\widetilde{Q}$ satisfies (\ref{8}) then for a fixed value of 
$s$ the matrix function $\widetilde{Q}$ is a solution to the $n$ waves problem.
The system (\ref{8}) admits the Lax representation
$$
 [L_i,L_j]=0, \ \ \
 L_j=\partial_{u^j}+I_j\partial_s+[I_j,\widetilde{Q}].
$$

The dressing method consist in the following procedure.
Let us consider the integral equation of the Marchenko type:
\beq
\label{9}
 K(s,s',u)=F(s,s',u)+\int_s^{\infty}K(s,q,u)F(q,s',u)dq,
\eeq 
where $F(s,s',u)$ is a matrix function such that it satisfies the equation
\beq
\label{10}
 \partial_{u^i}F+I_i\partial_sF+\partial_{s'}FI_i=0
\eeq 
and (\ref{9}) has a unique solution.
Then the function
\beq
\label{11}
\tilde{Q}(s,u) = K(s,s,u)
\eeq
satisfies (\ref{8}) and therefore for any fixed value of $s$ the function
$Q(u) = \widetilde{Q}(s,u)$
satisfies (\ref{7}). 
Moreover if the differential reduction
\beq
\label{12}
 \partial_{s'}F_{ij}(s,s',u)+\partial_{s}F_{ji}(s',s,u)=0,
\eeq 
holds (here $F_{ij}$ are the entries of the matrix $F$) then
$\widetilde{Q}$
also satisfies (\ref{5}).

The system consisting of (\ref{10}) and (\ref{12}) 
admits the following solution \cite{Z}:

\begin{itemize}
\item
{\sl Let $\Phi_{ij}(x,y), i<j$, be arbitrary $\frac{n(n-1)}{2}$
functions of two variables and 
$\Phi_{ii}(x,y)$ be $n$ arbitrary skew-symmetric functions:
$$
 \Phi_{ii}(x,y)=-\Phi_{ii}(y,x).
$$
We put
$$
 F_{ij}=\partial_s\Phi_{ij}(s-u^i,s'-u^j),\ i<j,
$$
$$
 F_{ji}=\partial_s\Phi_{ij}(s'-u^i,s-u^j),
$$
$$
 F_{ii}=\partial_s\Phi_{ii}(s-u^i,s'-u^i).
$$
The matrix function 
$F=(F_{ij})$ satisfies (\ref{10}) and (\ref{12}) and therefore
a solution $K$ to (\ref{9}) with such a matrix $F$ gives for any fixed value 
of $s$ the rotation coefficients of an orthogonal coordinate system:
$Q_{ij}(u) = K_{ij}(s,s,u)$, i.e. a solution to 
the system (\ref{4})--(\ref{5}).}
\end{itemize}

{\sc Remark 1.}
Since we have $\frac{n(n-1)}{2} + n = \frac{n(n+1)}{2}$ functional parameters,
i.e., $\Phi_{ij}, i \leq j$, and a general solution
depends on $\frac{n(n-1)}{2}$ functional parameters, this method gives 
equivalence classes of dressings
as it is explained in \cite{Z}.

\subsection{Krichever's method \cite{K}}
\label{subsec2.3}

Let $\Gamma$ be a smooth connected complex algebraic curve. We take three 
divisors on $\Gamma$:
$$
 P=P_1+\dots+P_n, \ \
D = \gamma_1+\dots+\gamma_{g+l-1}, \ \
 R=R_1+\dots+R_l,
$$
where $g$ is the genus of $\Gamma$, $P_i,\gamma_j,R_k\in\Gamma$.
We denote by 
$k_i^{-1}$ some local parameter on $\Gamma$ near $P_i$, $i=1,\dots,n$.

The Baker--Akhiezer function corresponding to the data
$$
 S=\{P,D,R\}
$$
is a function $\psi(u^1,\dots,u^n,z),\ z\in\Gamma$, 
meeting the following conditions:

1) $\psi\exp(-u^i k_i)$ is analytic near $P_i$, $i=1,\dots,n$;

2) $\psi$ is meromorphic on $\Gamma\backslash\{\cup P_i\}$ with poles at
$\gamma_j$, $j=1,\dots,g+l-1$;

3) $\psi(u,R_k)=1$, $k =1,\dots,l$.

For a generic divisor $D$ such a function exists and unique. Moreover it is 
expressed in terms of the theta function of $\Gamma$
\cite{K}. 

If the curve $\Gamma$ is not connected it is assumed that 
the restriction of the 
Baker--Akhiezer function onto every connected component meets
the conditions above.

We take an additional divisor 
$Q=Q_1+\dots+Q_n$ on $\Gamma$ such that $Q_i \in \Gamma \setminus 
\{P \cup D \cup R\}, i=1,\dots,n$, and denote by
$x^j$ the following function
$$
x^j(u^1,\dots,u^n)=\psi(u^1,\dots,u^n,Q_j), \ j =1,\dots,n.
$$
There is the following Krichever theorem
\cite{K}:

\begin{itemize}
\item
{\sl Let $\Gamma$ admit a holomorphic involution
$\sigma:\Gamma\rightarrow\Gamma$ 
such that 

1) this involution has exactly $2m, m \leq n$, fixed points 
which are just the points $P_1,\dots, P_n$ from $P$ and 
$2m-n$ points from $Q$;

2) $\sigma(Q)=Q$, i.e, non-fixed points from $Q$ are 
interchanged by the involution:
$$
\sigma(Q_k) = Q_{\sigma(k)}, \ \ \ k=1,\dots,n;
$$

3) $\sigma(k_i^{-1}) = -k_i^{-1}$ near $P_i$, $i=1,\dots,n$;

4) there exists a meromorphic differential
$\Omega$ on $\Gamma$ such that
its divisors of zeroes and poles are of the form
$$
(\Omega)_0= D + \sigma D +P, \ \ \ (\Omega)_{\infty}=R+\sigma R+Q.
$$

It is assumed that $\Gamma_0 = \Gamma/\sigma$ is a smooth algebraic curve.

Then, as it is easy to show, 
$\Omega$ is a pull-back of some meromorphic differential
$\Omega_0$ on $\Gamma_0$ and the following equalities hold:
$$
\sum_{k,l} \eta_{kl}\partial_{u^i}x^k\partial_{u^j} x^l
= \varepsilon^2_i h^2_i \delta_{ij}
$$
where
$$
h_i = \lim_{P \to P_i} \left(\psi e^{-u^i k_i}\right),
\ \ \ 
\eta_{kl} = \delta_{k,\sigma(l)} \res_{Q_k} \Omega_0,
\footnote{
Given a nonsingular
fixed point $Q_i$ of $\sigma$, there is a parameter 
$k$ near it such that $\sigma(k)=-k, k(Q_i)=0$. 
Therefore $\lambda = k^2$ is a local
parameter near $Q_i$ on $\Gamma/\sigma$ and we have 
$\Omega = \left(\frac{a}{k} + \dots\right)dk, 
\Omega_0 = \frac{1}{2} \left(\frac{a}{\lambda} + \dots\right) d\lambda$,
which implies
$\res_{Q_i} \Omega_0 = \frac{1}{2} \res_{Q_i}\Omega$.}
$$
and
$$
\Omega_0 = \frac{1}{2} \left(\varepsilon_i^2 \lambda_i + O(\lambda_i)\right) 
d\lambda_i, \ \lambda_i = k_i^{-2},
\ \mbox{ at $P_i$, $i=1,\dots,n$}.
$$
}
\end{itemize}

{\sc Remark 2.} 
This theorem stays true if instead of 1) we assume that 
the functions $\psi\exp(- f^i(u^i)k_i)$ are analytic near $P_i$
where $f^i$ are some functions of one variable, $i=1,\dots,n$.
Therefore we do not differ orthogonal coordinate systems which are obtained by 
coordinate changes of the form
$$
 u^i\rightarrow f^i(u^i).
$$

{\sc Remark 3.} 
Krichever's theorem gives a construction of these coordinates by using
the formalism of Baker--Akhiezer functions. In fact, it is clear from the proof
that the uniqueness of the Baker--Akhiezer function stays valid if we replace 
the condition $\psi(u,R_k) = 1, i=1,\dots,l$, by 
\beq
\label{baker1}
\psi(u,R_k) =d_k, \ \ k=1,\dots,l, 
\eeq
where all constants $d_k$ do not vanish. From that we deduce 
that we even can assume only that not all constants $d_k$ vanish:
\beq
\label{baker2}
|d_1|^2 + \dots + |d_l|^2 \neq 0,
\eeq
and under this assumption the main results of \cite{K} still hold. 

For distinguishing the cases when this Theorem give positive-definite 
metrics it needs to impose
some other conditions on the spectral data \cite{K}:

\begin{itemize}
\item
{\sl If there is an antiholomorphic involution $\tau: \Gamma \to \Gamma$
such that all fixed points of $\sigma$ are fixed by $\tau$ and
$$
\tau^\ast(\Omega) = \overline{\Omega}
$$
(for that it is enough to assume
that $\tau(k_i^{-1}) = \overline{k_i^{-1}}$
at $P_i, i=1,\dots,n$, and $\tau$ maps divisors $Q,R$, and $D$ into themselves:
$\tau(Q) = Q, \tau(R) = R, \tau(D) = D$, however these divisors do not 
necessarily consist of fixed points of $\tau$),
then the coefficients $H_i(u)$ are real valued for $u^1,\dots,u^n \in \R$.
}

\item
{\sl  $u^1,\dots,u^n$ are 
$n$-orthogonal coordinates 
in the flat $n$-space
with the metric $\eta_{kl} dx^k dx^l$.}
\footnote{It is easy to see that if there are points $Q_k$ and 
$Q_l, l =\sigma(k)$, 
which are interchanged by the involution then the metric $\eta$ is indefinite.}

\item
{\sl Provided that all points from $Q$ are fixed by the involution $\sigma$ and
\beq
\label{diff}
\res_{Q_1} \Omega_0 = \dots = \res_{Q_n} \Omega_0 = \eta^2_0 > 0,
\eeq
the functions $x^1(u^1,\dots,u^n),\dots,x^n(u^1,\dots,x^n)$ solve the immersion
problem for $n$-orthogonal coordinates $u^1,\dots,u^n$ and
$ds^2 = H_1^2 (du^1)^2 + \dots + H_n^2 (du^n)^2$ with
$$
H_i = \frac{\varepsilon_i h_i}{\eta_0}, \ \ \ i=1,\dots,n.
$$
}
\end{itemize}

The analog of Krichever's construction for discrete orthogonal systems was 
developed in \cite{AVK}.

Krichever's results allow us to assume that $n$-orthogonal coordinate systems
which are expressed in terms of elementary functions 
correspond to limiting cases when the spectral curve is singular.

\section{Coordinate systems corresponding to singular spectral curves}

Let $\Gamma$ be an algebraic complex curve with singularities.
Then there exists a morphism of a nonsingular algebraic curve
$\Gamma_\nm$:
$$
\pi: \Gamma_\nm \to \Gamma,
$$
such that

1) there is a finite set $S$ of points from $\Gamma_\nm$ 
with the equivalence relation $\sim$ on this set
such that $\pi$ maps $S$ exactly into the set 
$\Sing = \Sing
\Gamma$ formed by all singular points of $\Gamma$, 
and therewith the preimage of every point from 
$\Sing$ consists in a class of equivalent points;

2) the mapping $\pi: \Gamma_\nm \setminus S \to \Gamma \setminus
\Sing$ is a smooth one-to-one projection;

3) any regular mapping 
$F: X \to \Gamma$ of a nonsingular algebraic variety 
$X$ with an everywhere dense image $F(X) \subset \Gamma$
descends through $\Gamma_\nm$, i.e., $F = \pi G$ for some regular mapping
$G: X \to \Gamma_\nm$.

A mapping $\pi$ meeting these properties is called the normalization of
$\Gamma$ and is unique. The genus of
$\Gamma_\nm$ is called the geometric genus of
$\Gamma$ and is denoted by $p_g(\Gamma)$.

However another genus comes into the Riemann--Roch formula, i.e., 
the arithmetic genus $p_a(\Gamma)$ which is a sum of the geometric genus and 
some positive-valued contribution of singularities (the points from $\Sing$). 
For a nonsingular curve we have
$p_a = p_g$.

For example, let us consider the case of multiple points when on 
$\Gamma_\nm$ we choose $s$ families  
$D_1,\dots,D_s$ consisting of $r_1,\dots,r_s$ points all of which are 
pairwise different. Let $\Gamma$ is obtained by gluing together
points from each family. Then
$$
p_a(\Gamma) = p_g(\Gamma) + \sum_{i=1}^s (r_i-1).
$$
A meromorphic $1$-form $\omega$ on $\Gamma_\nm$ defines a regular differential
on  $\Gamma$ if for any point $P \in \Sing$ we have
$$
\sum_{\pi^{-1}(P)} \res \ (f\omega) = 0
$$
for any meromorphic function $f$, on $\Gamma_\nm$, 
which descends to a function on $\Gamma$, i.e., takes the same value
at points from each divisor $D_i$, and does not have poles in
$\pi^{-1}(P)$. Regular differentials may have poles in 
the preimages of singular sets.
It is easy to notice that the dimension of the space of regular differentials
equals  $p_a(\Gamma)$.

In the general case all these notions are exposed in
\cite{Serre} (for using it in the finite gap integration we gave some short
expositions in \cite{T1,T2}). We only remind the 
Riemann--Roch theorem for singular curves.

Let $L(D)$ be the space of meromorphic functions on 
$\Gamma$ with poles at the points from
$D = \sum n_P P$ of order less or equal that $n_P$,
and let $\Omega(D)$ be the space of regular differentials on
$\Gamma$  which has at every point $P \in \Sing$ a zero of order at least
$n_P$. The Riemann--Roch theorem reads
$$
\dim L(D) - \dim \Omega(D) = \deg D + 1 - p_a(\Gamma).
$$
For generic divisor $D$ with $\deg D > g$ we have $\dim
\Omega(D)=0$ and the Riemann--Roch theorem takes the form
$$
\dim L(D) = \deg D + 1 - p_a(\Gamma).
$$

\begin{theorem}
\label{th1}
Krichever's theorem (see \S \ref{subsec2.3}) holds for 
singular algebraic curves
provided that
$g$ is replaced by $p_a(\Gamma)$, i.e., by 
the arithmetic genus of $\Gamma$, and
the assumption that $\Gamma/\sigma$ is a nonsingular curve is replaced by 
the condition that $P_1,\dots,P_n$ and the poles of $\Omega$ are
nonsingular points.

Moreover we may assume that $\psi$ meets the conditions (\ref{baker1}) and
(\ref{baker2}) instead of $\psi(u,R_k)=1, k=1,\dots,l$. 
\end{theorem}

{\it Proof} 
of this theorem is basically the same as the original Krichever proof 
in \cite{K}. The uniqueness of $\psi$ is established
by using the general theory of Baker--Akhiezer functions. In the cases studied
in \S\S \, 3.1, 3.2 and 4 such a uniqueness is trivial since we are working 
with rational curves. 
The identity
$$
\sum_{k,l} \eta_{kl} \partial_{u^i} x^k \partial_{u^j} x^l -
\varepsilon_i^2 h^2_i \delta_{ij} = 0
$$
is equivalent to the identity
$$
\sum \res 
\left(\partial_{u^i} \psi(u,X) \partial_{u^j} \psi(u,\sigma(X))
\Omega \right) =0
$$
and is obtained from it by explicit calculations of the residues.

{\sc Remark 4 (main).} In the case when $\Gamma_\nm$ is 
a union of smooth rational curves, i.e., copies of $\C P^1$, the procedure of
constructing Baker--Akhiezer functions and orthogonal coordinates 
is very simple: it reduces to simple computations with elementary functions 
and does not go far than solving systems of linear equations.
However singular curves of algebraic genus $g$ are obtained via degeneration 
from smooth curves of the same genus. Therewith qualitative properties of
solutions, which correspond to these curves, of nonlinear equations are 
inherited, i.e., such solutions are rather complicated.

We restrict  ourselves by the most simple case when $\Gamma$ is a reducible 
curve consisting of components  $\Gamma_1,\dots,\Gamma_s$ isomorphic to 
$\C P^1$. These components may intersect each other at some points. 

A regular differential $\Omega$ on $\Gamma$ is defined by some differentials 
$\Omega_1,\dots,\Omega_s$ on
$\Gamma_1,\dots,\Gamma_s$ which may have poles at the intersections of 
components and moreover if $P$ is such an intersection point of the components
$\Gamma_{i_1},\dots,\Gamma_{i_r}$ then
$$
 \sum_{k=1}^r \res_P \ \Omega_{i_k}=0.
$$
The arithmetic genus $g_a$ is the dimension of the space of 
holomorphic regular differentials, i.e., differentials such that
$\Omega_j$ may have poles only at intersections of different components.

For different combinatorial configurations of 
rational components and intersection points Theorem \ref{th1}
is written in absolutely elementary form and 
the construction of orthogonal coordinates reduces
to solving some systems of linear equations.
It is simpler to demonstrate that by some examples which we expose below.
For simplicity we assume that on every component there is defined some 
complex parameter.

\subsection{$2$-orthogonal coordinate systems}

{\sc Example 1.}
Let $\Gamma$ consists of two copies of $\C P^1$, i.e., of $\Gamma_1$ and 
$\Gamma_2$, which intersect each other at two points:
$$
a \sim b, \ \ (-a) \sim (-b), \ \ \
 \{a,-a\} \subset \Gamma_1, \{b,-b\} \subset \Gamma_2
$$
(see Fig. 1). We have $p_a(\Gamma) = 1$.

The Baker--Akhiezer function takes the form
$$
 \psi_1(u^1,u^2,z_1)=e^{u^1z_1}\left(f_0(u^1,u^2)+
\frac{f_1(u^1,u^2)}{z_1-\alpha_1}+\dots
 +\frac{f_k(u^1,u^2)}{z_1-\alpha_{s_1}}\right),\ z_1\in\Gamma_1,
$$
$$
 \psi_2(u^1,u^2,z_2)=e^{u^2z_2}\left(g_0(u^1,u^2)+
\frac{g_1(u^1,u^2)}{z_2-\beta_1}+\dots
 +\frac{g_n(u^1,u^2)}{z_2-\beta_{s_2}}\right),\ z_2\in\Gamma_2.
$$
\beq
\label{gluing}
 \psi_1(a)=\psi_2(b),\ \ \psi_1(-a)=\psi_2(-b).
\eeq
It has two essential singularities at the points
$P_1=\infty\in\Gamma_1$ and $P_2=\infty\in\Gamma_2$.

\vskip15mm

\begin{picture}(170,100)(-100,-80)
\qbezier(-10,0)(-10,30)(40,30)
\qbezier(40,30)(90,30)(90,0)
\qbezier(90,0)(90,-30)(40,-30)
\qbezier(40,-30)(-10,-30)(-10,0)
\put(35,33){\shortstack{$\Gamma_1$}}

\qbezier(60,0)(60,30)(110,30)
\qbezier(110,30)(160,30)(160,0)
\qbezier(160,0)(160,-30)(110,-30)
\qbezier(110,-30)(60,-30)(60,0)
\put(105,33){\shortstack{$\Gamma_2$}}

\put(-10,0){\circle*{3}}
\put(90,0){\circle*{3}}
\put(60,0){\circle*{3}}
\put(160,0){\circle*{3}}
\put(75,24){\circle*{3}}
\put(75,-24){\circle*{3}}

\put(-25,0){\shortstack{$P_1$}}
\put(165,0){\shortstack{$Q_2$}}
\put(95,0){\shortstack{$Q_1$}}
\put(45,0){\shortstack{$P_2$}}
\put(63,20){\shortstack{\small{$a$}}}
\put(54,-26){\shortstack{\small{$-a$}}}
\put(85,18){\shortstack{\small{$b$}}}
\put(80,-26){\shortstack{\small{$-b$}}}
\put(60,-50){\shortstack{Fig. 1.}}
\end{picture}

The general normalization condition takes the form
\beq
\label{normal}
\psi_1(R_{1,i})= d_{1,i}, \ \
\psi_2(R_{2,j})= d_{2,j}
\eeq where
$R_{1,i}\in\Gamma_1, i=1,\dots,l_1$, and
$R_{2,j} \in \Gamma_2, j=1, \dots,l_2$.
We also have
$$
l = l_1 + l_2 = s_1 + s_2.
$$
Let
$$
\Omega_1=\frac{(z_1^2-\alpha_1^2)\dots (z_1^2-\alpha_{l_1}^2)}
{z_1(z_1^2-a^2)(z_1^2-R_{1,1}^2)\dots(z_1^2-R_{1,l_1}^2)}dz_1,
$$
$$
\Omega_2=\frac{(z_2^2-\beta_1^2)\dots (z_2^2-\beta_{l_2}^2)}
{z_2(z_2^2-b^2)(z_2^2-R_{2,1}^2)\dots(z_2^2-R_{2,l_2}^2)}dz_2.
$$
We put
$$
Q_1 = 0\in\Gamma_1,\ \ \ \ Q_2=0\in\Gamma_2.
$$
If the following equalities hold
$$
\res_a \Omega_1 = - \res_b \Omega_2,\ \
\res_{-a}\Omega_1 = -\res_{-b}\Omega_2,\ \
\res_{Q_1}\Omega_1 = \res_{Q_2}\Omega_2,
$$
then the differential $\Omega$ defined by 
$\Omega_1$ and $\Omega_2$ is regular, the condition (\ref{diff}) is satisfied 
and,
by Theorem \ref{th1}, the coordinates
$u^1$ and $u^2$ such that
$$
x^1(u) = \psi_1(u,0), \ \ x^2(u) = \psi_2(u,0)
$$
are orthogonal.

Let us consider the simplest case 
$l_1=0$ and $l_2 = 1$.

We have
$$
\psi_1=e^{u^1z_1}f_0(u^1,u^2),\ \
\psi_2=e^{u^2z_2}\left(g_0(u^1,u^2)+\frac{g_1(u^1,u^2)}{z_2-c}\right).
$$

\vskip25mm

\begin{picture}(170,100)(-100,-80)
\put(60,-30){\vector(0,4){130}}
\put(-15,25){\vector(3,0){185}}

\put(60,35){\circle{20}}
\put(60,45){\circle{40}}

\put(70,25){\circle{20}}
\put(80,25){\circle{40}}

\put(160,30){\shortstack{$x^1$}}
\put(65,92){\shortstack{$x^2$}}
\put(50,-50){\shortstack{Fig. 2.}}
\end{picture}

The gluing conditions at the intersection points and the 
normalization condition are
$$
\psi_1(a)=\psi_2(b),\ \ \psi_1(-a)=\psi_2(-b),\ \ \psi_2(r)=1, \ 
r = R \in \Gamma_2,
$$
which implies
$$
 \psi_1=e^{u^1z_1}\left(
  \frac{2b(c-r)e^{au^1+(b-r)u^2}}{(b+c)(b-r)e^{2bu^2}-(b+r)(b-c)e^{2au^1}}
  \right),
$$
$$
 \psi_2=e^{u^2z_2}\left(
 \frac{e^{-ru^2}((b-c)e^{2au^1}+(b+c)e^{2bu^2})(c-r)}
 {(b+c)(b-r)e^{2bu^2}-(b-c)(b+r)e^{2au^1}}+\right.
$$
$$
 \left. \frac{1}{z_2-c}\frac{(b^2-c^2)(r-c)e^{-ru^2}(e^{2au^1}-e^{2bu^2})}
 {(b+c)(r-b)e^{2bu^2}+(b-c)(b+r)e^{2au^1}}\right).
$$

The differential $\Omega$ is defined by
the differentials
$$
\Omega_1 = -\frac{1}{z_1(z_1^2-a^2)}dz_1,\
\Omega_2 = - \frac{(z_2^2-c^2)}
{z_2(z_2^2-b^2)(z_2^2-r^2)}dz_2.
$$
We have the following  regularity condition for $\Omega$:
$$
\res_{a}\Omega_1=\res_{-a}\Omega_1 = -\frac{1}{2a^2}= 
- \res_{b}\Omega_2= -\res_{-b}\Omega_2= 
\frac{(b^2-c^2)}{2b^2(b^2-r^2)},
$$
and the condition (\ref{diff}) takes the form
$$
\res_{Q_1}\Omega_1= \frac{1}{a^2}=
\res_{Q_2}\Omega_2 = \frac{c^2}{r^2 b^2}
$$
which implies
$$
 a=\frac{br}{c},\ r=\frac{b}{\sqrt{2 + \frac{b^2}{c^2}}}.
$$

After the substitution
$u^1=\log y^1,\ u^2=\log y^2$
the formulas for the coordinates 
are written as
$$
x^1=\frac{-2b(r-c)}{(c-b)(b+r)}(y^2)^{-r}\left(\frac{\frac{(y^2)^{b}}{(y^1)^a}}
{1+\frac{(b+c)(b-r)}{(c-b)(b+r)} \frac{(y^2)^{2b}}{(y^1)^{2a}}}\right),
$$
$$
x^2=\frac{b(c-r)}{c(b+r)}(y^2)^{-r}
\left(\frac{1+\frac{(b+c)}{(c-b)} \frac{(y^2)^{2b}}{(y^1)^{2a}}}
{{1+\frac{(b+c)(b-r)}{(c-b)(b+r)} \frac{(y^2)^{2b}}{(y^1)^{2a}}}}\right)
$$
and by straightforward computations we obtain
$$
 (x^1)^2+\left(x^2-(y^2)^{-r}\frac{b(c-r)}{c(b^2-r^2)}\right)^2=
 (y^2)^{-2r}\frac{b^2(c-r)^2}{c^2(b^2-r^2)^2}.
$$
Therefore the coordinate lines $y^2=\const$, i.e., $u^2=\const$, are 
the circles centered on the $x^2$ axis. For $b = \pm 1$ these circles touch 
the $x^1$ axis and another family of coordinate lines consist of 
circles centered at the $x^1$ axis and touching the $x^2$ axis
(see Fig. 2).

{\sc Example 2.} 
Let  $\Gamma$ be the same as in Example 1 however
all essential singularities lie in one copy of $\C P^1$ and the divisor 
$Q$ lies in another copy (see Fig. 3)
$$
P_1=\infty, \ P_2=0\in\Gamma_1,\ \ Q_1=\infty, \ Q_2=0 \in \Gamma_2.
$$
We define the Baker--Akhiezer function as follows:
$$
 \psi_1(u,z_1)=e^{u^1z_1+\frac{u^2}{z_1}}\left(f_0(u)+
\frac{f_1(u)}{z_1-\alpha_1}+\dots
 +\frac{f_k(u)}{z_1-\alpha_{s_1}}\right), \ z_1\in\Gamma_1;
$$
$$ 
\psi_2(u,z_2)=\left(g_0(u)+\frac{g_1(u)}{z_2-\beta_1}+\dots
 +\frac{g_n(u)}{z_2-\beta_{s_2}}\right), \ z_2\in\Gamma_2.
$$

\vskip2cm

\begin{picture}(170,100)(-100,-80)
\qbezier(-10,0)(-10,30)(40,30)
\qbezier(40,30)(90,30)(90,0)
\qbezier(90,0)(90,-30)(40,-30)
\qbezier(40,-30)(-10,-30)(-10,0)

\put(35,33){\shortstack{$\Gamma_1$}}

\qbezier(60,0)(60,30)(110,30)
\qbezier(110,30)(160,30)(160,0)
\qbezier(160,0)(160,-30)(110,-30)
\qbezier(110,-30)(60,-30)(60,0)

\put(105,33){\shortstack{$\Gamma_2$}}

\put(-10,0){\circle*{3}}
\put(90,0){\circle*{3}}
\put(60,0){\circle*{3}}
\put(160,0){\circle*{3}}
\put(75,24){\circle*{3}}
\put(75,-24){\circle*{3}}

\put(-25,0){\shortstack{$P_1$}}
\put(165,0){\shortstack{$Q_2$}}
\put(95,0){\shortstack{$P_2$}}
\put(45,0){\shortstack{$Q_1$}}
\put(63,20){\shortstack{\small{$a$}}}
\put(54,-26){\shortstack{\small{$-a$}}}
\put(85,18){\shortstack{\small{$b$}}}
\put(80,-26){\shortstack{\small{$-b$}}}
\put(60,-50){\shortstack{Fig. 3.}}

\end{picture}

The gluing and normalization conditions have the same forms
(\ref{gluing}) and (\ref{normal}).
Let
$$
\Omega_1=\frac{z_1(z_1^2-\alpha_1^2)\dots (z_1^2-\alpha_{l_1-1}^2)}
{(z_1^2-a^2)(z_1^2-R_{1,1}^2)\dots(z_1^2-R_{1,l_1}^2)}dz_1,
$$
$$
\Omega_2=\frac{(z_2^2-\beta_1^2)\dots (z_2^2-\beta_{l_2+1}^2)}
{z_2(z_2^2-b^2)(z_2^2-R_{2,1}^2)\dots(z_2^2-R_{2,l_2}^2)}dz_2.
$$

By Theorem \ref{th1}, if
$$
 \res_a\Omega_1=- \res_b\Omega_2,\
 \res_{-a}\Omega_1=- \res_{-b}\Omega_2,\
 \res_{Q_1}\Omega_2 = \res_{Q_2}\Omega_2,
$$
then we have
$$
 \partial_{u^1}x^1\partial_{u^2}x^1+\partial_{u^1}x^2\partial_{u^2}x^2=0.
$$

Let us consider the simplest case:
$l_1=s_2 = 1, l_2=s_1=0$, $r=R \in \Gamma_1$, $d_{1,1}=1$.
We have
$$
 \psi_1=e^{u^1z_1+\frac{u^2}{z_1}}f(u),\ \
\psi_2=\left(g_0(u)+\frac{g_1(u)}{z_2-c}\right),
$$
$$
\Omega_1=\frac{z_1}{(z_1^2-a^2)(z_1^2-r^2)}dz_1, \ \
\Omega_2= - \frac{(z_2^2-c^2)}
{z_2(z_2^2-b^2)}dz_2.
$$
We have
$$
\res_{a}\Omega_1= \res_{-a}\Omega_1=\frac{1}{2(a^2-r^2)}, \ \ 
\res_{b}\Omega_2= \res_{-b}\Omega_2=\frac{(b^2-c^2)}{2b^2},
$$
$$
\res_{Q_1}\Omega_2=1, \ \ \ \res_{Q_2}\Omega_2=-\frac{b^2}{c^2},
$$
and the regularity condition for $\Omega$ and (\ref{diff}) are satisfied 
exactly when
$$
b = \pm ic, \ \ \ a^2 - r^2 = -\frac{1}{2}.
$$
For a particular solution 
$b=i,\ c=-1,\ a=\frac{i}{2},\ r=\frac{1}{2}$,
the immersion formulas take the form
$$
x^1=e^{-\frac{u^1}{2}-2u^2}\left(\cos\left({\frac{u^1}{2}-2u^2}\right)+
\sin\left({\frac{u^1}{2}-2u^2}\right)
\right),
$$
$$
x^2=e^{-\frac{u^1}{2}-2u^2}\left(\cos\left({\frac{u^1}{2}-2u^2}\right)-
\sin\left({\frac{u^1}{2}-2u^2}\right)
\right).
$$
By the substitution
$$
 y^1=\frac{u^1}{2},\ y^2=2u^2
$$
we obtain
$$
x^1=e^{-y^1-y^2}(\cos(y^1-y^2)+
\sin(y^1-y^2)),
$$
$$
x^2=e^{-y^1-y^2}(\cos(y^1-y^2)-
\sin(y^1-y^2)).
$$
Therewith the ``lines'' 
$y^1+y^2=\const$ correspond to circles centered at the origin $x=0$,
and the ``lines'' $y^1-y^2=\const$ define in the $x$-space rays 
drawing from the origin.

\subsection{$3$-orthogonal coordinate systems}

{\sc Example 3.} Let $\Gamma$ consist of three components
$\Gamma_1$, $\Gamma_2$ and $\Gamma_3$ which are copies of $\C P^1$ and
have four intersection points as it is shown on Fig. 4:

\vskip1cm

\begin{picture}(170,100)(-100,-80)
\qbezier(-60,0)(-60,30)(-10,30)
\qbezier(-10,30)(40,30)(40,0)
\qbezier(40,0)(40,-30)(-10,-30)
\qbezier(-10,-30)(-60,-30)(-60,0)
\put(-15,33){\shortstack{\footnotesize$\Gamma_1$}}

\qbezier(10,0)(10,30)(60,30)
\qbezier(60,30)(110,30)(110,0)
\qbezier(110,0)(110,-30)(60,-30)
\qbezier(60,-30)(10,-30)(10,0)
\put(55,33){\shortstack{\footnotesize$\Gamma_2$}}

\qbezier(80,0)(80,30)(130,30)
\qbezier(130,30)(180,30)(180,0)
\qbezier(180,0)(180,-30)(130,-30)
\qbezier(130,-30)(80,-30)(80,0)
\put(125,33){\shortstack{\footnotesize$\Gamma_3$}}

\put(-60,0){\circle*{3}}
\put(40,0){\circle*{3}}
\put(40,0){\circle*{3}}
\put(80,0){\circle*{3}}
\put(10,0){\circle*{3}}

\put(110,0){\circle*{3}}
\put(25,24){\circle*{3}}
\put(25,-24){\circle*{3}}
\put(180,0){\circle*{3}}
\put(95,24){\circle*{3}}
\put(95,-24){\circle*{3}}

\put(-75,0){\shortstack{\footnotesize$P_1$}}
\put(115,0){\shortstack{\footnotesize$Q_1$}}
\put(45,0){\shortstack{\footnotesize$P_2$}}
\put(-5,0){\shortstack{\footnotesize$P_3$}}

\put(185,0){\shortstack{\footnotesize$Q_3$}}
\put(64,0){\shortstack{\footnotesize$Q_2$}}

\put(14,20){\shortstack{\small{$a$}}}
\put(5,-26){\shortstack{\small{$-a$}}}
\put(35,18){\shortstack{\small{$b$}}}
\put(30,-26){\shortstack{\small{$-b$}}}

\put(84,20){\shortstack{\small{$c$}}}
\put(74,-26){\shortstack{\small{$-c$}}}
\put(105,18){\shortstack{\small{$d$}}}
\put(100,-26){\shortstack{\small{$-d$}}}

\put(50,-50){\shortstack{Fig. 4}}

\end{picture}

$$
\pm a \sim \pm b, \ \ \pm c  \sim \pm d, \ \ \pm a \in\Gamma_1, \ \pm b, \pm c 
\in \Gamma_2, \ \ \pm d \in \Gamma_3.
$$

Let us put 
$$
P_1 = \infty \in \Gamma_1, \ \ P_2 = 0 \in \Gamma_1, \ \
P_3 = \infty \in \Gamma_2, 
$$
$$
Q_1 = 0 \in \Gamma_2, \ \
Q_2 = \infty \in \Gamma_3, \ \ Q_3 = 0 \in \Gamma_3, 
$$
$$
l=1, \ \ r = R \in \Gamma_1, \ \ \psi_1(r) =1.
$$

We have
$$
 \psi_1(u,z_1)=e^{u^1z_1+\frac{u^2}{z_1}}f(u),\ \ z_1\in\Gamma_1,
$$
$$
 \psi_2(u,z_2)=e^{u^3z_2}
\left(g_0(u)+\frac{g_1(u)}{z_2-\beta}\right),\ \ z_2\in\Gamma_2,
$$
$$
 \psi_3(u,z_3)=h_0(u)+\frac{h_1(u)}{z_3-\gamma},\ \ z_3\in\Gamma_3.
$$
with
$$
\psi_1(\pm a)=\psi_2(\pm b),\ \psi_2(\pm c)=\psi_3(\pm d),\ \
\psi_1(r) =1.
$$

Let us take $\Omega$ defined by the differentials
$$
\Omega_1= \frac{z_1dz_1}
{(z_1^2-a^2)(r^2 - z_1^2)}, \ \
\Omega_2 = \frac{(\beta^2 - z_2^2)dz_2}
{z_2(z_2^2-b^2)(z_2^2-c^2)}, \ \
\Omega_3 = \frac{(\gamma^2- z_3^2)dz_3}
{z_3(z_3^2-d^2)}.
$$

The regularity condition for $\Omega$ and (\ref{diff}) are satisfied 
for
$$
a^2 = -\frac{1}{12}, \ b = -\frac{1}{3}, \ c=d=i, \ \beta= bc, \ \gamma=-1
$$
in which case we have
$$
 x^1=\sqrt{2}e^{-\frac{u^1}{2}-2u^2}\cos
\left(\frac{1}{12}(3\pi+2\sqrt{3}(u^1-2(6u^2+u^3)))\right),
$$
$$
 x^2=\sqrt{2}e^{-\frac{u^1}{2}-2u^2}
 \left(\cos\left(\frac{u^1-2(6u^2+u^3)}{2\sqrt{3}}\right)
\sin\left(\frac{\pi}{4}+u^3\right)+
\right.
$$
$$
\left.
\sin\left(\frac{u^1-2(6u^2+u^3)}{2\sqrt{3}}\right)
\cos\left(\frac{\pi}{12}+u^3\right)\right),
$$
$$
 x^3=\sqrt{2}e^{-\frac{u^1}{2}-2u^2}
 \left(\cos\left(\frac{u^1-2(6u^2+u^3)}{2\sqrt{3}}\right)
\cos\left(\frac{\pi}{4}+u^3\right)-
\right.
$$
$$
\left.
\sin\left(\frac{u^1-2(6u^2+u^3)}{2\sqrt{3}}\right)
\sin\left(\frac{\pi}{12}+u^3\right)\right).
$$
It is straightforwardly checked that
$$
 (x^1)^2+(x^2)^2+(x^3)^2=3e^{-u^1-4u^2},
$$
$$
x^1-\left(\frac{1-\sqrt{3}}{2}x^2+\frac{1+\sqrt{3}}{2}x^3\right)\cos u^3-
\left(\frac{1+\sqrt{3}}{2}x^2+ \frac{\sqrt{3}-1}{2}x^3\right)\sin u^3=0,
$$
$$
 2(x^1)^2-(x^2)^2-(x^3)^2+((x^1)^2+(x^2)^2+(x^3)^2)
\sin\left(\frac{u^1-2(6u^2+u^3)}{\sqrt{3}}\right)=0.
$$
Therefore the ``planes'' $u^3=\const$ are planes passing the point $x=0$,
the ``planes'' $u^1+4u^2=\const$ are spheres centered at $x=0$, and
the ``planes'' $u^1-2(6u^2+u^3)=\const$ are cones centered at $x=0$.

\section{The classical coordinate systems}

{\bf The Euclidean coordinates.}
Let $\Gamma$ be a disjoint union of $n$ copies $\Gamma_1,\dots,\Gamma_n$ 
of $\C P^1$. We put
$$
P_j=\infty, \ \ Q_j=0, \ \ R_j=-1\in\Gamma_j, \ \ \psi_j(R_j)=1, \ \ \ 
j=1,\dots,n.
$$
Then we have the differential $\Omega$ defined by the 
differentials
$$
\Omega_j=\frac{dz_j}{z_j (z_j^2-1)}, \ \ \ j=1,\dots,n,
$$
on the components of $\Gamma$. The Baker--Akhiezer function $\psi$ 
is equal to
$$
\psi_j=e^{u^jz_j}f_j(u^j), \ \ \ j=1,\dots,n,
$$
and we obtain the Euclidean coordinates in $\R^n$ (see Remark 2):
$$
x^j = e^{u^j}.
$$

\medskip

{\bf The polar coordinates.}
Let $\Gamma$ consists of five irreducible components  $\Gamma_1,\dots,
\Gamma_5$ which intersect as it is
shown on Fig. 5:
$$
\{0\in\Gamma_1\}\sim\{0\in\Gamma_2\},\ \ 
\{a\in\Gamma_2\}\sim\{b_1\in\Gamma_3\},
\ \ \{-a\in\Gamma_2\}\sim\{b_2\in\Gamma_4\},
$$
$$
 \{c_1\in\Gamma_3\}\sim\{d\in\Gamma_5\}, \ \
 \{c_2\in\Gamma_4\}\sim\{-d\in\Gamma_5\}.
$$
We define an involution $\sigma$ on $\Gamma$ as follows:

a) on $\Gamma_1,\Gamma_2$ and $\Gamma_3$ it has the form
$$
\sigma(z_j)=-z_j;
$$

b) $\Gamma_3$ and $\Gamma_4$ are interchanged by $\sigma$ and the points
$b_1,c_1,\infty\in\Gamma_3$ are mapped into the points
$b_2,c_2,\infty\in\Gamma_4$ respectively. It is easy to check that
$$
\sigma(z_3)=\frac{b_2-c_2}{b_1-c_1}z_3+\frac{b_1c_2-b_2c_1}{b_1-c_1},
$$
$$
\sigma(z_4)=\frac{b_1-c_1}{b_2-c_2}z_4+\frac{b_2c_1-b_1c_2}{b_2-c_2}.
$$
We put
$$
\beta_1=\frac{b_2c_1-b_1c_2}{b_2-c_2},\
\beta_2=\frac{b_1c_2-b_2c_1}{b_1-c_1}.
$$
Then $0\in\Gamma_3$ is mapped by $\sigma$ into
$\beta_2\in\Gamma_4$, and $0\in\Gamma_4$ is mapped into
$\beta_1\in\Gamma_3$.

The divisors $P =P_1+P_2$ and $Q=Q_1+Q_2$ are as follows
$$
 P_1=\infty\in\Gamma_1,\ \ 
P_2=\infty\in\Gamma_2,\ \ 
Q_1 = 0 \in \Gamma_5,\ \ 
Q_2 = \infty \in \Gamma_5.
$$
We also take for $D = \gamma_1 + \gamma_2 + \gamma_3$ the divisor 
$$
\gamma_1 = 0 \in \Gamma_3,\ \ \gamma_2 =  0 \in \Gamma_4, \ \ 
\gamma_3 = \alpha \in \Gamma_5.
$$
Since $\deg D = g+l-1 = 3$ and $g=p_a(\Gamma) =1$, we have $l=3$.
We put
$$
R_1 = -1 \in \Gamma_1, \ \  
R_2 = \infty \in \Gamma_4, \ \ 
R_3 = \infty \in \Gamma_5.
$$

Then the Baker--Akhiezer function takes the form
$$
 \psi_1(u,z_1)=e^{u^1 z_1}f_1(u),\ \
 \psi_2(u,z_2)=e^{u^2 z_2}f_2(u), 
$$
$$
 \psi_3(u,z_3)=\frac{f_3(u)}{z_3}+\widehat{f}_3(u), \ \
 \psi_4(u,z_4)=\frac{f_4(u)}{z_4}+\widehat{f}_4(u), 
$$
$$ 
 \psi_5(u,z_5)=f_5(u)+\frac{\widehat{f}_5(u)}{z_5-\alpha}.
$$

\vskip30mm

\begin{picture}(85,50)(-150,-80)
\qbezier(-100,0)(-100,20)(-70,20)
\qbezier(-70,20)(-40,20)(-40,0)
\qbezier(-40,0)(-40,-20)(-70,-20)
\qbezier(-70,-20)(-100,-20)(-100,0)

\qbezier(-50,0)(-50,40)(-20,40)
\qbezier(-20,40)(10,40)(10,0)
\qbezier(10,0)(10,-40)(-20,-40)
\qbezier(-20,-40)(-41,-37)(-46,-23)

\qbezier(7,39)(14,44)(30,45)
\qbezier(30,45)(47,44)(54,39)
\qbezier(60,25)(60,5)(30,5)
\qbezier(30,5)(0,5)(0,25)

\qbezier(0,-25)(0,-5)(30,-5)
\qbezier(30,-5)(60,-5)(60,-25)
\qbezier(54,-39)(47,-44)(30,-45)
\qbezier(30,-45)(14,-44)(7,-40)

\qbezier(50,0)(50,40)(80,40)
\qbezier(80,40)(110,40)(110,0)
\qbezier(110,0)(110,-40)(80,-40)
\qbezier(80,-40)(50,-40)(50,0)

\put(-100,0){\circle*{3}}

\put(-49,15){\circle*{3}}
\put(-43,15){\shortstack{\footnotesize$0$}}
\put(-57,7){\shortstack{\footnotesize$0$}}

\put(10,0){\circle*{3}}
\put(14,-3){\shortstack{\footnotesize$P_2$}}

\put(50,0){\circle*{3}}
\put(37,-3){\shortstack{\footnotesize$Q_1$}}

\put(110,0){\circle*{3}}
\put(113,-4){\shortstack{\footnotesize$Q_2$}}

\put(10,8){\circle*{3}}
\put(2,3){\shortstack{\footnotesize$a$}}
\put(14,10){\shortstack{\footnotesize$b_1$}}

\put(10,-9){\circle*{3}}
\put(-6,-9){\shortstack{\footnotesize$-a$}}
\put(14,-17){\shortstack{\footnotesize$b_2$}}

\put(50,8){\circle*{3}}
\put(41,10){\shortstack{\footnotesize$c_1$}}
\put(54,1){\shortstack{\footnotesize$d$}}

\put(50,-9){\circle*{3}}
\put(41,-17){\shortstack{\footnotesize$c_2$}}
\put(53,-10){\shortstack{\footnotesize$-d$}}

\put(-113,0){\shortstack{\footnotesize$P_1$}}
\put(-72,25){\shortstack{$\Gamma_1$}}
\put(-25,45){\shortstack{$\Gamma_2$}}

\put(25,50){\shortstack{$\Gamma_3$}}
\put(25,-40){\shortstack{$\Gamma_4$}}
\put(78,45){\shortstack{$\Gamma_5$}}

\put(15,-60){\shortstack{Fig. 5}}

\end{picture}

Moreover we have
$$
 \psi_1(u,0)=\psi_2(u,0),\ \ \psi_2(u,a)=\psi_3(u,b_1),\ \
 \psi_2(u,-a)=\psi_4(u,b_2),
$$
$$
 \psi_3(u,c_1)=\psi_5(u,d),\ \
 \psi_4(u,c_2)=\psi_5(u,-d).
$$
We take the following normalization condition:
$$
 \psi_1(u,-1)=1,\ \ \psi_3(u,\infty)=0,\
 \psi_4(u,\infty)=0.
$$

This implies
$$
f_1=e^{u^1},\ \ f_2=e^{u^1},\ \
f_3 = e^{u^1+au^2},\ \ 
\widehat{f}_3=0,\ \ 
f_4 = e^{u^1-au^2},\ \ 
\widehat{f}_4=0,
$$
$$
 f_5 = 
\frac{e^{u^1-au^2}(b_1c_2e^{2au^2}(d-\alpha)+b_2c_1(d+\alpha))}{2c_1c_2d},
$$
$$
 \widehat{f}_5=\frac{e^{u^1-au^2}(-b_2c_1+e^{2au^2}
 b_1c_2)(d^2-\alpha^2)}{2c_1c_2d}.
$$

Is is checked by straightforward computations that for
$$
a = i, \ \ b_1 = \bar{b}_2 = \frac{i}{2}, \ \ 
c_1 = \bar{c}_2 = \frac{i-1}{2}, \ \ d=-i\alpha
$$
the differential $\Omega$ defined by the differentials
$$
\Omega_1 = -\frac{dz_1}{z_1 (z_1^2-1)},\
\Omega_2 = -\frac{dz_2}{z_2(z_2^2-a^2)},\
\Omega_3 = -\frac{z_3(z_3-\beta_1)dz_3}{(z_3-b_1)(z_3-c_1)},
$$
$$
\Omega_4 = -\frac{z_4(z_4-\beta_2)dz_3}{(z_4-b_2)(z_4-c_2)},\
\Omega_5 = -\frac{(z_5^2-\alpha^2)dz_5}{z_5(z_5^2-d^2)}
$$
is regular on $\Gamma$, satisfies (\ref{diff}), and
$$
 x^1 = \psi_5(Q_1) = r \cos \varphi, \ \ 
 x^2 = \psi_5(Q_2) =  r \sin \varphi
$$
where $r = e^{u^1}$ and $\varphi = u^2$.

{\sc Remark 5.} We see that the values of $\alpha$ and $d$ are not precisely 
determined and we only have the relation $\alpha = id$. Therefore as in 
the case of the dressing method this construction also corresponds 
an equivalence  class of spectral data to the same metric 
(see Remark 1.)

\medskip

{\bf The cylindrical coordinates.}
We take for $\Gamma$ a disjoint union of the curve $\widehat{\Gamma}$
from the previous example
(the polar coordinates) and a copy $\Gamma_6$ of $\C P^1$. 
All the data concerning $\widehat{\Gamma}$ are also the same as in 
the previous example. On $\Gamma_6$ we put $Q_3 = 0, P_3 = \infty, R_4 = -1$,
 and $\psi(R_4) = 1$. 
Then  we have
$\psi_6(u^3)=e^{u^3(z_6+1)}$ and
$$
x^1 = \psi_5(Q_1) = r \cos \varphi, \ \ 
x^2 = \psi_5(Q_2) = r \sin \varphi, \ \
x^3 = \psi_6(Q_3) = z
$$
where $r = e^{u^1}, \varphi= u^2$, and $z = u^3$.

\medskip

{\bf The spherical coordinates in ${\mathbb R}^3$.} 
The curve $\Gamma$ consists of $9$ irreducible components which 
intersect as it is shown on Fig. 6:

\vskip20mm

\begin{picture}(170,100)(-90,-80)

\qbezier(-85,0)(-85,20)(-55,20)
\qbezier(-55,20)(-25,20)(-25,0)
\qbezier(-25,0)(-25,-20)(-55,-20)
\qbezier(-55,-20)(-85,-20)(-85,0)

\put(-85,0){\circle*{3}}
\put(-96,-4){\shortstack{\scriptsize$P_1$}}
\put(-65,25){\shortstack{\footnotesize$\Gamma_1$}}
\put(-42,18){\circle*{3}}
\put(-50,10){\shortstack{\footnotesize$0$}}
\put(-37,20){\shortstack{\footnotesize$0$}}

\qbezier(-45,0)(-45,40)(-15,40)
\qbezier(-15,40)(15,40)(15,0)
\qbezier(15,0)(15,-40)(-15,-40)
\qbezier(-15,-40)(-36,-37)(-41,-23)

\put(-15,45){\shortstack{\footnotesize$\Gamma_2$}}
\put(15,0){\circle*{3}}
\put(15,7){\circle*{3}}
\put(15,-7){\circle*{3}}
\put(18,-3){\shortstack{\scriptsize$P_2$}}
\put(7,2){\shortstack{\footnotesize$a$}}
\put(0,-7){\shortstack{\footnotesize$-a$}}
\put(18,8){\shortstack{\footnotesize$b_1$}}
\put(19,-15){\shortstack{\footnotesize$b_2$}}

\qbezier(7,39)(14,44)(30,45)
\qbezier(30,45)(47,44)(54,39)
\qbezier(60,25)(60,5)(30,5)
\qbezier(30,5)(0,5)(0,25)
\put(27,47){\shortstack{\footnotesize$\Gamma_3$}}
\put(45,7){\circle*{3}}
\put(36,8){\shortstack{\footnotesize$c_1$}}

\qbezier(0,-25)(0,-5)(30,-5)
\qbezier(30,-5)(60,-5)(60,-25)
\qbezier(54,-39)(47,-44)(30,-45)
\qbezier(30,-45)(14,-44)(7,-40)
\put(27,-40){\shortstack{\footnotesize$\Gamma_4$}}
\put(45,-7){\circle*{3}}
\put(36,-15){\shortstack{\footnotesize$c_2$}}

\qbezier(45,0)(45,40)(75,40)
\qbezier(75,40)(105,40)(105,0)
\qbezier(105,0)(103,-29)(97,-31)
\qbezier(75,-40)(45,-40)(45,0)
\put(45,0){\circle*{3}}
\put(33,-3){\shortstack{\scriptsize$Q_1$}}
\put(48,1){\shortstack{\footnotesize$d$}}
\put(47,-8){\shortstack{\footnotesize$-d$}}

\put(72,43){\shortstack{\footnotesize$\Gamma_5$}}
\put(90,36){\circle*{3}}
\put(80,31){\shortstack{\footnotesize$0$}}
\put(97,31){\shortstack{\footnotesize$0$}}

\qbezier(75,0)(75,40)(105,40) \qbezier(105,40)(135,40)(135,0)
\qbezier(135,0)(135,-40)(105,-40) \qbezier(105,-40)(75,-40)(75,0)
\put(102,43){\shortstack{\footnotesize$\Gamma_6$}}
\put(135,0){\circle*{3}} \put(135,7){\circle*{3}}
\put(135,-7){\circle*{3}}
\put(139,-3){\shortstack{\scriptsize$P_3$}}
\put(126,2){\shortstack{\footnotesize$a$}}
\put(120,-7){\shortstack{\footnotesize$-a$}}
\put(138,9){\shortstack{\footnotesize$b_1$}}
\put(138,-15){\shortstack{\footnotesize$b_2$}}

\qbezier(127,39)(134,44)(150,45) \qbezier(150,45)(167,44)(174,39)
\qbezier(180,25)(180,5)(150,5) \qbezier(150,5)(120,5)(120,25)
\put(150,47){\shortstack{\footnotesize$\Gamma_7$}}
\put(165,7){\circle*{3}}
\put(155,8){\shortstack{\footnotesize$c_1$}}

\qbezier(120,-25)(120,-5)(150,-5) \qbezier(150,-5)(180,-5)(180,-25)
\qbezier(174,-39)(167,-44)(150,-45)
\qbezier(150,-45)(134,-44)(127,-40)
\put(150,-40){\shortstack{\footnotesize$\Gamma_8$}}
\put(165,-7){\circle*{3}}
\put(155,-14){\shortstack{\footnotesize$c_2$}}

\qbezier(165,0)(165,40)(195,40)
\qbezier(195,40)(225,40)(225,0)
\qbezier(225,0)(225,-40)(195,-40)
\qbezier(195,-40)(165,-40)(165,0)

\put(165,0){\circle*{3}} \put(152,-3){\shortstack{\scriptsize$Q_2$}}
\put(169,1){\shortstack{\footnotesize$d$}}
\put(168,-8){\shortstack{\footnotesize$-d$}}
\put(225,0){\circle*{3}} \put(228,-3){\shortstack{\scriptsize$Q_3$}}
\put(192,43){\shortstack{\footnotesize$\Gamma_9$}}

\put(75,-60){\shortstack{Fig. 6}}

\end{picture}

$$
\{0\in\Gamma_1\}\sim\{0\in\Gamma_2\},\ \ 
\{a\in\Gamma_2\}\sim\{b_1\in\Gamma_3\},\ \ 
\{-a\in\Gamma_2\}\sim\{b_2\in\Gamma_4\},
$$
$$
\{c_1\in\Gamma_3\}\sim\{d\in\Gamma_5\},\ \ 
\{c_2\in\Gamma_4\}\sim\{-d\in\Gamma_5\},\ \ 
\{0\in\Gamma_5\}\sim\{0\in\Gamma_6\},\ \  
$$
$$
\{a\in\Gamma_6\}\sim\{b_1\in\Gamma_7\},\ \ 
\{-a\in\Gamma_6\}\sim\{b_2\in\Gamma_8\},\ \
$$
$$
\{c_1\in\Gamma_7\}\sim\{d\in\Gamma_9\},\ \ 
\{c_2\in\Gamma_8\}\sim\{-d\in\Gamma_9\}
$$
where, for simplicity, 
we denote by the same symbol points on different components if the
coordinates of these points are equal to each other (for, instance, $a$ 
on $\Gamma_2$ and $\Gamma_6$).

We take
$$
Q_1=\infty\in \Gamma_5 ,\ \ Q_2=\infty \in \Gamma_9,\ \ Q_3=0\in\Gamma_9,
$$
$$
P_1 = \infty \in \Gamma_1, \ \ P_2 = \infty \in \Gamma_2, \ \ 
P_3 = \infty \in \Gamma_6,
$$
and choose the divisor $D$ as follows
$$
\gamma_1 = 0 \in \Gamma_3, \ \ \gamma_2 = 0 \in \Gamma_4, \ \
\gamma_3 = \alpha \in \Gamma_5,
$$
$$
\gamma_4 = 0 \in \Gamma_7, \ \ \gamma_5 = 0 \in \Gamma_8, \ \
\gamma_6 = \alpha \in \Gamma_9.
$$
We have $p_a(\gamma)=2$, $\deg D = 6$, and therefore $l = 5$.
Let us put
$$
R_1 = -1 \in \Gamma_1, \ \ R_2 = \infty \in \Gamma_3, \ \ 
R_3 = \infty \in \Gamma_4, R_4 = \infty \in \Gamma_7, \ \ 
R_5 = \infty \in \Gamma_8.
$$

The Baker--Akhiezer function is written as
$$
\psi_1=e^{u^1z_1}f_1(u),\ \
\psi_2=e^{u^2z_2}f_2(u), \ \
\psi_3=\frac{f_3(u)}{z_3}+\widehat{f}_3(u),
$$
$$
\psi_4=\frac{f_4(u)}{z_4}+\widehat{f}_4(u),\ \
\psi_5=f_5(u)+\frac{\widehat{f}_5(u)}{(z_5-\alpha)},\ \
\psi_6=e^{u^3z_6}f_6(u),
$$
$$
\psi_7=\frac{f_7(u)}{z_7}+\widehat{f}_7(u),\ \
\psi_8=\frac{f_8(u)}{z_8}+\widehat{f}_8(u),\ \
\psi_9=f_9(u)+\frac{\widehat{f}_9(u)}{z_9-\alpha}.
$$
We have the gluing conditions (for brevity,we skip the $u$-variables):
$$
 \psi_1(0)=\psi_2(0),\ \ \psi_2(a)=\psi_3(b_1),\ \
 \psi_2(-a)=\psi_4(b_2),\ \
 \psi_3(c_1)=\psi_5(d),\
$$
$$
 \psi_4(c_2)=\psi_5(-d),\ \ \psi_5(0)=\psi_6(0),\ \
 \psi_6(a)=\psi_7(b_1),\ \
 \psi_6(-a)=\psi_8(b_2),\
$$
$$
 \psi_7(c_1)=\psi_9(d),\ \ \
 \psi_8(c_2)=\psi_9(-d),
$$
and the normalization condition is taken as follows:
$$
\psi_1(u,-1)=1, \ \psi_3(u,\infty)=0,\ \psi_4(u,\infty)=0, \
\psi_7(u,\infty)=0,\ \psi_8(u,\infty)=0.
$$

Let $a,b_1,b_2,c_1,c_2,d$ take the same values as in the case of polar 
coordinates and then the regular form $\Omega$ is also 
constructed as in the this case.
By straightforward computations, we obtain
$$
x^1 = \psi_5(Q_1) = r \sin \varphi,
$$
$$
x^2 = \psi_9(Q_2) = r \cos \varphi \sin \theta, \ \
x^3 = \psi_9(Q_3) = r \cos \varphi \cos \theta
$$
where $r = e^{u^1}, \varphi = u^2$, and $\theta = u^3$.

\medskip

{\bf The spherical coordinates in ${\mathbb R}^n$.}
Let $\Gamma^{(n-1)}$ be the spectral curve and $\psi^{(n-1)}$
be the Baker--Akhiezer function 
for the $(n-1)$-dimensional spherical coordinates.
Then the spectral curve $\Gamma^{(n)}$ for 
the $n$-dimensional spherical coordinates
is the union of $\Gamma^{(n-1)}$ and the curve exposed on Fig. 7. Therewith
these curves intersect at the points $0\in\Gamma_{4n-7} \subset \Gamma^{(n-1)}$
and $0\in\Gamma_{4n-6}$ (we remark that the number of irreducible components 
of $\Gamma^{(k)}$ equals $4k-3$). The Baker--Akhiezer function coincides with
$\psi^{(n-1)}$ on $\Gamma^{(n-1)}$ and on the additional components
it is defined as follows:
$$
P_n=\infty\in\Gamma_{4n-6},\ \ 
Q_{n-1}=\infty,\ \ 
Q_n=0\in\Gamma_{4n-3},
$$
and on the additional components the Baker--Akhiezer function is defined 
as follows:
$$
\psi_{4n-6}=e^{u^n z_{4n-6}}f_{4n-6}(u),\ \ 
\psi_{4n-5}=\frac{f_{4n-5}(u)}{z_{4n-5}},
$$
$$
 \psi_{4n-4}=\frac{f_{4n-4}(u)}{z_{4n-4}}, \ \ \
 \psi_{4n-3}=f_{4n-3}(u)+\frac{\widehat{f}_{4n-3}(u)}{z_{4n-3}-\alpha}.
$$

\vskip30mm

\begin{picture}(85,50)(-150,-80)

\qbezier(-50,0)(-50,40)(-20,40)
\qbezier(-20,40)(10,40)(10,0)
\qbezier(10,0)(10,-40)(-20,-40)
\qbezier(-20,-40)(-50,-40)(-50,0)
\put(-50,0){\circle*{3}}
\put(-46,-3){\shortstack{\footnotesize$0$}}

\qbezier(7,39)(14,44)(30,45)
\qbezier(30,45)(47,44)(54,39)
\qbezier(60,25)(60,5)(30,5)
\qbezier(30,5)(0,5)(0,25)

\qbezier(0,-25)(0,-5)(30,-5)
\qbezier(30,-5)(60,-5)(60,-25)
\qbezier(54,-39)(47,-44)(30,-45)
\qbezier(30,-45)(14,-44)(7,-40)

\qbezier(50,0)(50,40)(80,40)
\qbezier(80,40)(110,40)(110,0)
\qbezier(110,0)(110,-40)(80,-40)
\qbezier(80,-40)(50,-40)(50,0)

\put(10,0){\circle*{3}}
\put(14,-3){\shortstack{\scriptsize$P_n$}}

\put(50,0){\circle*{3}}
\put(28,-3){\shortstack{\scriptsize$Q_{n-1}$}}

\put(110,0){\circle*{3}}
\put(113,-3){\shortstack{\scriptsize$Q_n$}}

\put(10,8){\circle*{3}} \put(3,2){\shortstack{\footnotesize$a$}}
\put(14,10){\shortstack{\footnotesize$b_1$}}

\put(10,-9){\circle*{3}} \put(-5,-8){\shortstack{\footnotesize$-a$}}
\put(14,-17){\shortstack{\footnotesize$b_2$}}

\put(50,8){\circle*{3}} \put(41,10){\shortstack{\footnotesize$c_1$}}
\put(54,1){\shortstack{\footnotesize$d$}}

\put(50,-9){\circle*{3}}
\put(41,-17){\shortstack{\footnotesize$c_2$}}
\put(53,-10){\shortstack{\footnotesize$-d$}}

\put(-27,45){\scriptsize{$\Gamma_{4n-6}$}}

\put(23,50){\scriptsize{$\Gamma_{4n-5}$}}
\put(23,-40){\scriptsize{$\Gamma_{4n-4}$}}
\put(76,45){\scriptsize{$\Gamma_{4n-3}$}}

\put(15,-60){\shortstack{Fig. 7}}

\end{picture}

\end{document}